\documentclass[11pt]{article}%
\usepackage{amsmath}%
\usepackage{amsfonts}%
\usepackage{amssymb}%
\usepackage{graphicx}

\newtheorem{theorem}{Theorem}

\newtheorem{algorithm}[theorem]{Algorithm}

\newtheorem{definition}[theorem]{Definition}
\newtheorem{example}[theorem]{Example}

\newtheorem{lemma}[theorem]{Lemma}
\newtheorem{notation}[theorem]{Notation}

\newtheorem{proposition}[theorem]{Proposition}
\newtheorem{remark}[theorem]{Remark}

\newenvironment{proof}[1][Proof]{\noindent\textbf{#1.} }{\ \rule{0.5em}{0.5em}}
\begin{document}
\title{Parametrization of the regular equivalences of the
canonical controller and its applications}
\author{A. Agung Julius\footnote{Dept. Electrical and Systems
Engineering, Univ. Pennsylvania, 200 S 33rd Street, Philadelphia,
PA 19104, USA. Email:agung@seas.upenn.edu}%
\; Jan Willem Polderman\footnote{Dept. Applied Mathematics, Univ.
Twente, PO Box 217, Enschede 7500AE, The Netherlands.
Email:j.w.polderman@math.utwente.nl}%
\; and Arjan van der Schaft\footnote{Institute for Mathematics and
Computer Science, Univ. Groningen, PO Box 800, Groningen 9700AV,
The Netherlands. Email:A.J.van.der.Schaft@math.rug.nl}}

\maketitle
\begin{abstract}
We study control problems for linear systems in the behavioral
framework. Our focus is a class of regular controllers that are
equivalent to the canonical controller. The canonical controller
is a particular controller that is guaranteed to be a solution
whenever a solution exists.  However, it has been shown that in
most cases, the canonical controller is not regular. The main
result of the paper is a parametrization of all regular
controllers that are equivalent to the canonical controller. The
parametrization is then used to solve two control problems. The
first problem is related to designing a regular controller that
uses as few control channels as possible. The second problem is to
design a regular controller that satisfies a predefined
input-output partitioning constraint. In both problems, based on
the parametrization, we present algorithms that does the
controller design.
\end{abstract}

\section{Introduction}

In this paper, we discuss control problems for linear differential
systems in the behavioral approach. The behavior of a system is
the set of trajectories that are compatible with the laws that
describe the system. In the continuous time case, the behavior is
the set of solutions of the differential equations that describe
the system. In the discrete time case, it is the set of solutions
of some difference equations.

\begin{figure}
[ptb]
\begin{center}
\includegraphics[
height=3.15cm, width=6.0803cm
]%
{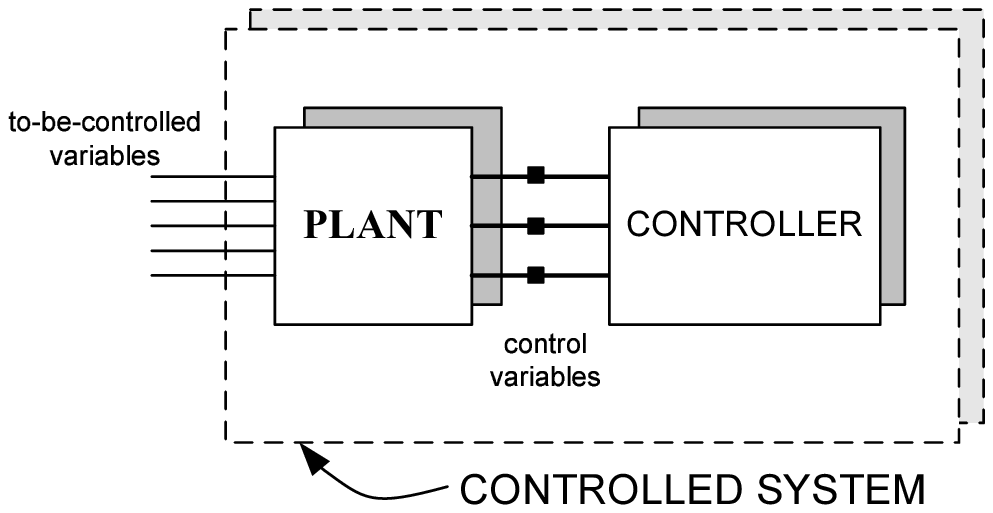}%
\caption{Control in the behavioral approach.}%
\label{f1}%
\end{center}
\end{figure}

Standard control problems in the behavioral approach to systems
theory can be formulated as follows
\cite{Willems97,Trentelman02,Belur03}. Given is a plant to be
controlled that has two kinds of variables: \emph{to-be-controlled
variables} and \emph{control variables}. A \emph{controller} is a
device that is attached to (or an algorithm that acts on) the
control variables and restricts their behavior. This restriction
is imposed on the plant via the control variables, such that it
(indirectly) affects the behavior of the to-be-controlled
variables (see Figure \ref{f1}). The resulting behavior is
called the \emph{controlled system}.%

As part of the control problem, one is given a
\emph{specification}, which is expressed in terms of the
to-be-controlled variables. The objective of the control problem
is to make the controlled system satisfy the specification. If
there exists a controller such that this objective is satisfied,
we say that the specification is \emph{implementable}.

In \cite{Julius02,Schaft03}, a particular controller design,
called the \emph{canonical controller }was introduced. This design
has the nice property that it implements the desired specification
if and only if the specification is implementable. However,
analysis on the regularity of the canonical controller reveals
that it is \emph{maximally irregular} \cite{Julius04b}. Regularity
is a desirable property for the interconnection
\cite{Willems97,Trentelman02}, which we will explain in Section
\ref{sec2}. In this paper, we show that there exist regular
controllers that are equivalent to the canonical controller, and
we provide a parametrization of all such controllers. This
parametrization is then used to solve two control problems:

\begin{enumerate}
\item The problem of \textbf{control with minimal interaction} \cite{Julius05}%
. This problem is about designing a regular controller that
interacts with the plant with as few control variables as
possible. The motivation behind this problem is as follows.
Consider a situation where the plant and the controller are
separated by a large physical distance. We need a communication
link between the plant and the controller to establish the
interconnection. It is therefore favorable to have as few control
variables as possible, so that the amount of communication
links/channels can be minimized.

\item The problem of \textbf{control with I/O partitioning
constraint}. This problem is about designing a regular controller
that respects a constraint on the \emph{a priori} partitioning of
the control variables into input and output variables.
\end{enumerate}

The results in this paper are presented in the form of continuous
time systems. However, they also hold for discrete time systems,
as we replace the differential operator with the discrete time lag
operator.

\section{Background material\label{sec2}}

For linear differential systems, the plant is typically described
as a set of linear differential equations that relate the
variables. Throughout this paper, we denote the control variables
as $\mathbf{c}$ and the to-be-controlled variables as
$\mathbf{w}$. The dimensions of $\mathbf{c}$ and $\mathbf{w}$ are
denoted as $\mathtt{c}$ and $\mathtt{w}$ respectively. A
behavioral model of the plant system that captures the relevant
relation between $\mathbf{w}$ and $\mathbf{c}$ is called the
\emph{full plant behavior}, and is denoted by
$\mathcal{P}_{\text{full}}.$ The full plant
behavior can be compactly represented as follows.%
\begin{equation}
\left[
\begin{array}
[c]{cc}%
R\left(  \frac{d}{dt}\right)  & M\left(  \frac{d}{dt}\right)
\end{array}
\right]  \left[
\begin{array}
[c]{c}%
w\\
c
\end{array}
\right]  =0, \label{1}%
\end{equation}
where $R$ and $M$ are polynomial matrices with appropriate
dimensions. We denote the class of polynomial matrices with
indeterminate $\xi$, $g$ rows, and $q$ columns over the real field
as $\mathbb{R}^{g\times q}[\xi]$.

The representation of the behavior in the form of (\ref{1}) is
called a \emph{kernel representation}, the reason being that the
behavior is simply the kernel of a linear differential operator.
Kernel representations of a given behavior are not unique. For
example, if $U(d/dt)$ is a linear differential operator, whose
kernel consists of only the zero trajectory, then the kernel of
$U(d/dt)\circ R(d/dt)$ is the same as that of $R(d/dt)$. Square
polynomial
matrices $U(\xi)$ such that%
\begin{equation}
\ker U\left(  \frac{d}{dt}\right)  =\{0\},
\end{equation}
are called \emph{unimodular matrices}$.$ It can be proven that the
inverse of $U(\xi)$ is also a polynomial matrix. A stronger result
that relates unimodular matrices and equivalent kernel
representations is that two kernel representations $R_{1}(d/dt)$
and $R_{2}(d/dt)$ with the same number of rows are equivalent if
and only if there is a unimodular matrix $U(\xi)$ such that
$R_{1}(\xi)=U(\xi)R_{2}(\xi)$.

Although the kernel representation of a behavior $\mathfrak{B}$ is
not unique, there is a unique integer $\mathtt{p}(\mathfrak{B)},$
which is the minimum number of rows a kernel representation of
$\mathfrak{B}$ can have. This number is also the row rank of any
kernel representation of the behavior. A kernel representation
with the minimum number of rows (i.e. equal to its row rank) is
called a \emph{minimal kernel representation}. The number $\mathtt{p}%
(\mathfrak{B})$ is called the number of outputs of $\mathfrak{B}$.

Suppose that a behavior $\mathfrak{B}$ is given by%
\begin{equation}
\mathfrak{B}:=\left\{  w~|~R\left(  \frac{d}{dt}\right)
w=0\right\}  ,
\label{13}%
\end{equation}
where $R$ is full row rank and has $\mathtt{p}(\mathfrak{B})$
rows. We can partition the variables in $\mathbf{w}$ into
$\mathbf{w}_{1}$ and
$\mathbf{w}_{2}$ such that (\ref{13}) becomes%
\begin{equation}
\mathfrak{B}:=\left\{  (w_{1},w_{2})~|~R_{1}\left(
\frac{d}{dt}\right) w_{1}+R_{2}\left(  \frac{d}{dt}\right)
w_{2}=0\right\}  ,
\end{equation}
where $R_{1}$ is a square full row rank polynomial matrix. Such a
partition is called an \emph{input-output} partition where
$\mathbf{w}_{1}$ is the output and $\mathbf{w}_{2}$ is the input
to the system. Notice that the number of outputs of $\mathfrak{B}$
is $\mathtt{p}(\mathfrak{B})$.

In this paper, we restrict our attention to infinitely
differentiable functions. Thus, the full plant behavior consists
of all signal pairs $(w,c)$ that are strong solutions to the
kernel representation (\ref{1})
\cite{Willems98}.%
\begin{equation}
\mathcal{P}_{\text{full}}:=\left\{  (w,c)\in\mathfrak{C}^{\infty}%
(\mathbb{R},\mathbb{R}^{\mathtt{w}+\mathtt{c}})~|~R\left(  \frac{d}%
{dt}\right)  w+M\left(  \frac{d}{dt}\right)  c=0\right\}  .
\end{equation}
If we eliminate the control variables from the full behavior, we
obtain the so
called \emph{manifest behavior}, which is denoted by $\mathcal{P}.$ Thus,%
\begin{equation}
\mathcal{P}:=\{w\in\mathfrak{C}^{\infty}(\mathbb{R},\mathbb{R}^{\mathtt{w}%
})~|~\exists~c\in\mathfrak{C}^{\infty}(\mathbb{R},\mathbb{R}^{\mathtt{c}%
})\text{ such that }(w,c)\in\mathcal{P}_{\text{full}}\}.
\end{equation}
If we rewrite the kernel representation (\ref{1}) as%
\begin{equation}
\left[
\begin{array}
[c]{cc}%
\tilde{R}_{1}\left(  \frac{d}{dt}\right)  & \tilde{M}_{1}\left(  \frac{d}%
{dt}\right) \\
\tilde{R}_{2}\left(  \frac{d}{dt}\right)  & 0
\end{array}
\right]  \left[
\begin{array}
[c]{c}%
w\\
c
\end{array}
\right]  =0,
\end{equation}
where $\tilde{M}_{1}$ and $\tilde{R}_{2}$ are full row rank
matrices, then the manifest behavior $\mathcal{P}$ is the kernel
of $\tilde{R}_{2}\left( \frac{d}{dt}\right)  $ (cf.
\cite{Willems98} Chapter 6).

A controller $\mathcal{C}$ is a behavior containing all signals
$c$ allowed by
the controller:%
\begin{equation}
\mathcal{C}:=\left\{  c\in\mathfrak{C}^{\infty}(\mathbb{R},\mathbb{R}%
^{\mathtt{c}})~|~C\left(  \frac{d}{dt}\right)  c=0\right\}  .
\end{equation}
The \emph{controlled behavior} is then defined as%
\begin{equation}
\mathcal{K}:=\{w\in\mathfrak{C}^{\infty}(\mathbb{R},\mathbb{R}^{\mathtt{w}%
})~|~\exists~c\in\mathfrak{C}^{\infty}(\mathbb{R},\mathbb{R}^{\mathtt{c}%
})\text{ such that }(w,c)\in\mathcal{P}_{\text{full}}\text{ and }%
c\in\mathcal{C\}}.
\end{equation}
The controlled behavior $\mathcal{K}$ is obtained by eliminating
the control
variables from the following kernel representation.%
\begin{align}
R\left(  \frac{d}{dt}\right)  w+M\left(  \frac{d}{dt}\right)  c  &
=0,\nonumber\\
C\left(  \frac{d}{dt}\right)  c  &  =0.
\end{align}
The specification $\mathcal{S}$ is given by the following kernel
representation%
\begin{equation}
S\left(  \frac{d}{dt}\right)  w=0. \label{4}%
\end{equation}
The objective of the control problem is to find a controller
$\mathcal{C}$ such that $\mathcal{K}=\mathcal{S}$. If such
controller exists, then $\mathcal{S}$ is said to be implementable
and the controller $\mathcal{C}$ is said to implement
$\mathcal{S}$.

Clearly, the implementability of a specification $\mathcal{S}$ is
a property that depends on the specification itself as well as the
plant. The following result is proven in
\cite{Willems99a,Willems02a}.

\begin{theorem}
[Willems' lemma]Given $\mathcal{P}_{\text{full}}$ as a kernel
representation
of (\ref{1}). A specification $\mathcal{S}$ is implementable if and only if%
\begin{equation}
\mathcal{N}\subseteq\mathcal{S}\subseteq\mathcal{P}, \label{2}%
\end{equation}
where $\mathcal{N}\in\mathfrak{L}^{\mathtt{w}}$ is the
\emph{hidden behavior}
defined by%
\[
\mathcal{N}:=\{w\in\mathfrak{C}^{\infty}(\mathbb{R},\mathbb{R}^{\mathtt{w}%
})~|~(w,0)\in\mathcal{P}_{\text{full}}\}.
\]

\end{theorem}

Quite often, in addition to requiring that the controller
implements the desired specification, we also require that the
controller possesses a certain property with respect to the plant.
A property that has been quite extensively studied is the so
called \emph{regularity}
\cite{Polderman00,Trentelman02,Julius03,Willems03a}. A controller%
\begin{equation}
\mathcal{C}=\left\{  c\in\mathfrak{C}^{\infty}(\mathbb{R},\mathbb{R}%
^{\mathtt{c}})~|~C\left(  \frac{d}{dt}\right)  c=0\right\}  , \label{3}%
\end{equation}
where $C$ is full row rank, to be regular if%
\begin{equation}
\text{rank}\left[
\begin{array}
[c]{cc}%
R & M\\
0 & C
\end{array}
\right]  =\text{rank}\left[
\begin{array}
[c]{cc}%
R & M
\end{array}
\right]  +\text{rank }C. \label{5}%
\end{equation}
It can be proven that nonregular interconnections affect the
autonomous part of the systems \cite{Willems97}, which, in many
cases would be undesirable or unrealistic.

\begin{remark}
Although the characterization of regular controllers suggests that
regularity is a representation dependent property, it is actually
not. Notice that (\ref{5}) is equivalent to saying that the number
of outputs of the controlled system is the the sum of those of the
plant and the controller. The number of outputs of a system, as
discussed earlier in this section, is a representation independent
quantity. The interested readers are referred to
\cite{Trentelman02,Julius03,Willems03a} for more discussion on the
behavioral interpretation of regularity.
\end{remark}

If the specification $\mathcal{S}$ is such that there exists a
regular controller $\mathcal{C}$ that implements it, then
$\mathcal{S}$ is said to be regularly implementable. Necessary and
sufficient conditions for regular implementability were derived in
\cite{Trentelman02}:

\begin{theorem}
\label{t1}Given the full plant behavior
$\mathcal{P}_{\text{full}}.$ A specification $\mathcal{S}$ is
regularly implementable if and only
if\newline1) it is implementable, i.e., $\mathcal{N}\subseteq\mathcal{S}%
\subseteq\mathcal{P}$ and\newline2) $\mathcal{S}+\mathcal{P}^{\text{ctr}%
}=\mathcal{P}.$\newline The symbol $\mathcal{P}^{\text{ctr}}$
denotes the controllable part of the manifest behavior
$\mathcal{P}.$
\end{theorem}

\section{The canonical controller and its regular equivalences}

In this section, we review the idea of canonical controller and
its properties. Given a full plant behavior
$\mathcal{P}_{\mathrm{full}}$ and a specification $\mathcal{S}$,
the behavior of the canonical controller
$\mathcal{C}_{\mathrm{can}}$ is defined as%
\begin{equation}
\mathcal{C}_{\mathrm{can}}:=\{c\in\mathfrak{C}^{\infty}(\mathbb{R}%
,\mathbb{R}^{\mathtt{c}})~|~\exists~w\in\mathfrak{C}^{\infty}(\mathbb{R}%
,\mathbb{R}^{\mathtt{w}})\text{ such that }(w,c)\in\mathcal{P}%
_{\text{\textrm{full}}}\text{ and }w\in\mathcal{S\}}.
\end{equation}
A kernel representation of the canonical controller can be
obtained by
eliminating $\mathbf{w}$ from the following kernel representation%
\begin{equation}
\left[
\begin{array}
[c]{cc}%
R\left(  \frac{d}{dt}\right)  & M\left(  \frac{d}{dt}\right) \\
S\left(  \frac{d}{dt}\right)  & 0
\end{array}
\right]  \left[
\begin{array}
[c]{c}%
w\\
c
\end{array}
\right]  =0.
\end{equation}

The canonical controller has the following property.

\begin{theorem}
(cf. \cite{Schaft03}) The canonical controller
$\mathcal{C}_{\mathrm{can}}$ implements the specification
$\mathcal{S}$ if and only if $\mathcal{S}$ is implementable.
\end{theorem}

We define the \emph{control manifest behavior} of the plant,
$\mathcal{P}_{c}$
as%
\begin{equation}
\mathcal{P}_{c}:=\{c\in\mathfrak{C}^{\infty}(\mathbb{R},\mathbb{R}%
^{\mathtt{c}})~|~\exists~w\in\mathfrak{C}^{\infty}(\mathbb{R},\mathbb{R}%
^{\mathtt{w}})\text{ such that
}(w,c)\in\mathcal{P}_{\text{\textrm{full}}}\}.
\end{equation}
A kernel representation of $\mathcal{P}_{c}$ can be obtained by
eliminating
$\mathbf{w}$ from the kernel representation of $\mathcal{P}%
_{\text{\textrm{full}}}$. The canonical controller has the
property of being \emph{least restrictive} in the following sense.

\begin{proposition}
\label{p1}(cf. \cite{Schaft03}) Assume that the specification
$\mathcal{S}$ is implementable. For any controller $\mathcal{C}$
that implements $\mathcal{S}$,
we have that%
\begin{equation}
\left(  \mathcal{C}\cap\mathcal{P}_{c}\right)  \subseteq\left(  \mathcal{C}%
_{\mathrm{can}}\cap\mathcal{P}_{c}\right)  .
\end{equation}
Thus, any trajectory of the control variables of the plant allowed
by $\mathcal{C}$ is also allowed by $\mathcal{C}_{\mathrm{can}}$.
\end{proposition}

Another important property of the canonical controller that is
relevant to our discussion in this paper, is that it is
\emph{maximally irregular}, in the following sense.

\begin{theorem}
(cf. \cite{Julius04b}) Assume that the specification $\mathcal{S}$
is implementable. The canonical controller
$\mathcal{C}_{\mathrm{can}}$ is regular if and only if every
controller that implements $\mathcal{S}$ is regular.
\end{theorem}

Although the canonical controller is maximally irregular, there
are regular controllers that are equivalent to it. By equivalent
controllers, we mean the controllers that allow the same set of
$\mathbf{c}$ trajectories of the plant as the canonical controller
does. The class of such controllers is defined as follows.

\begin{definition}
The class of regular controllers that are equivalent to the
canonical controller is denoted as
$\mathfrak{C}_{\mathrm{can}}^{\mathrm{reg}}$, and is
defined as%
\begin{equation}
\mathfrak{C}_{\mathrm{can}}^{\mathrm{reg}}:=\{\mathcal{C}~|~\mathcal{C}\text{
is regular and }\left(  \mathcal{C}\cap\mathcal{P}_{c}\right)
=\left( \mathcal{C}_{\mathrm{can}}\cap\mathcal{P}_{c}\right)  \}.
\end{equation}

\end{definition}

The following theorem provides a necessary and sufficient
conditions for the nonemptyness of the class
$\mathfrak{C}_{\mathrm{can}}^{\mathrm{reg}}$.

\begin{theorem}
\label{t2}The class $\mathfrak{C}_{\mathrm{can}}^{\mathrm{reg}}$
is nonempty if and only if the specification $\mathcal{S}$ is
regularly implementable.
\end{theorem}

\begin{proof}
The (only if) part of the theorem is obvious. We shall prove the
(if) part. Suppose that $\mathcal{S}$ is regularly implementable.
There exists a regular controller that implements $\mathcal{S}$.
We denote this controller as $\mathcal{C}$. By definition, we have
that\newline(a) For all $w\in \mathcal{S}$, there exists a
$c\in\mathcal{C}$ such that $(w,c)\in
\mathcal{P}_{\text{\textrm{full}}}$.\newline(b) For all
$c\in\mathcal{C}$,
$(w,c)\in\mathcal{P}_{\text{\textrm{full}}}$ implies $w\in\mathcal{S}%
$.\newline

Define another controller%
\begin{equation}
\mathcal{C}^{\prime}:=\mathcal{C}+\mathcal{C}_{\mathrm{can}}.
\end{equation}
We shall prove that $\mathcal{C}^{\prime}\in\mathfrak{C}_{\mathrm{can}%
}^{\mathrm{reg}}$, that is\newline(a') $\mathcal{C}^{\prime}$ is
regular.\newline(b') $\mathcal{C}^{\prime}\cap\mathcal{P}_{c}=\mathcal{C}%
_{\mathrm{can}}\cap\mathcal{P}_{c}$.\newline The statement (a')
follows from the fact that
$\mathcal{C}\subset\mathcal{C}^{\prime}$ and the regularity of
$\mathcal{C}$. To prove (b'), first we show that
$\mathcal{C}^{\prime}$ implements $\mathcal{S}.$ From here, (b')
follows from the fact that
$\mathcal{C}_{\mathrm{can}}\subseteq\mathcal{C}^{\prime}$ and the
property of $\mathcal{C}_{\mathrm{can}}$ being the least
restrictive controller (see Proposition \ref{p1}).

Showing that $\mathcal{C}^{\prime}$ implements $\mathcal{S}$ means
showing that\newline(a\textquotedblright) For all
$w\in\mathcal{S}$, there exists a
$c^{\prime}\in\mathcal{C}^{\prime}$ such that $(w,c^{\prime})\in
\mathcal{P}_{\text{\textrm{full}}}$.\newline(b\textquotedblright)
For all
$c^{\prime}\in\mathcal{C}^{\prime}$, $(w,c^{\prime})\in\mathcal{P}%
_{\text{\textrm{full}}}$ implies $w\in\mathcal{S}$.\newline
Statement
(a\textquotedblright) follows from (a) and the fact that $\mathcal{C}%
\subset\mathcal{C}^{\prime}$. To show that (b\textquotedblright)
holds, notice that any $c^{\prime}\in\mathcal{C}^{\prime}$ can be
written as
$c+c_{\text{\textrm{can}}}$ with $c\in\mathcal{C}$ and $c_{\text{\textrm{can}%
}}\in\mathcal{C}_{\text{\textrm{can}}}$. Also notice that for all
$c_{\text{\textrm{can}}}\in\mathcal{C}_{\text{\textrm{can}}}$,
there exists a
$w_{\text{\textrm{can}}}\in\mathcal{S}$ such that $(w_{\text{\textrm{can}}%
},c_{\text{\textrm{can}}})\in\mathcal{P}_{\text{\textrm{full}}}$. Thus,%
\begin{gather*}
(w,c^{\prime})\in\mathcal{P}_{\text{\textrm{full}}}\Rightarrow
(w-w_{\text{\textrm{can}}}+w_{\text{\textrm{can}}},c+c_{\text{\textrm{can}}%
})\in\mathcal{P}_{\text{\textrm{full}}}\\
\overset{\text{\textrm{linearity}}}{\Rightarrow}((w-w_{\text{\textrm{can}}%
}),c)\in\mathcal{P}_{\text{\textrm{full}}}\\
\overset{\text{\textrm{(b)}}}{\Rightarrow}(w-w_{\text{\textrm{can}}}%
)\in\mathcal{S}\\
\overset{\text{\textrm{linearity}}}{\Rightarrow}w\in\mathcal{S}\text{.}%
\end{gather*}

\end{proof}

The proof of Theorem \ref{t2} also implies the following important
property of $\mathfrak{C}_{\mathrm{can}}^{\mathrm{reg}}$.

\begin{theorem}
\label{t3}Given a control problem with a regularly implementable
specification $\mathcal{S}$. If $\mathcal{C}$ is a regular
controller that implements
$\mathcal{S}$, then there exists a regular controller $\mathcal{C}^{\prime}%
\in\mathfrak{C}_{\mathrm{can}}^{\mathrm{reg}}$ that implements
$\mathcal{S}$ and $\mathcal{C}\subseteq\mathcal{C}^{\prime}$.
\end{theorem}

One of the main results presented in this paper is the
parametrization of all controllers in
$\mathfrak{C}_{\mathrm{can}}^{\mathrm{reg}}$. Before we can obtain
the parametrization, we need the following lemma.

\begin{lemma}
\label{l1}Let a plant $\mathcal{P}$ be given as the kernel of a
full row rank $R\left(  \frac{d}{dt}\right)  $ and a regular
controller $\mathcal{C}$ be given as the kernel of a full row rank
$C\left(  \frac{d}{dt}\right)  .$
Denote the full interconnection%
\[
\mathcal{K}:=\mathcal{P}\cap\mathcal{C}.
\]
Let $\mathfrak{C}_{\mathcal{K}}$ denote the set of all controllers
(not necessarily regular ones) that\newline(i) have at most as
many outputs as $\mathcal{C}$ and\newline(ii) also implement
$\mathcal{K}$ when interconnected with $\mathcal{P}$.\newline A
controller $\mathcal{C}^{\prime}\in \mathfrak{C}_{\mathcal{K}}$ if
and only if its kernel representation can be written as $VR+C$ for
some matrix $V.$ Moreover, every controller in
$\mathcal{C}^{\prime}\in\mathfrak{C}_{\mathcal{K}}$ has the
following properties.\newline(a) $\mathcal{C}^{\prime}$ is
regular.\newline(b) $\mathcal{C}^{\prime}$ has exactly as many
outputs as $\mathcal{C}$.
\end{lemma}

\begin{proof}
(if) Suppose that a controller $\mathcal{C}^{\prime}$ is the
kernel of
$(VR+C)$, then $\mathcal{P}\cap\mathcal{C}^{\prime}$ is given by the kernel of%
\begin{equation}
\left[
\begin{array}
[c]{c}%
R\\
VR+C
\end{array}
\right]  =\left[
\begin{array}
[c]{cc}%
I & 0\\
V & I
\end{array}
\right]  \left[
\begin{array}
[c]{c}%
R\\
C
\end{array}
\right]  .
\end{equation}
This shows that $\mathcal{P\cap C}^{\prime}=\mathcal{P}\cap\mathcal{C}%
=\mathcal{K}$. Moreover, since $\mathcal{C}$ is a regular
controller, it
follows that $(VR+C)$ is a minimal kernel representation of $\mathcal{C}%
^{\prime}.$ Thus, properties (a) and (b) are verified.

(only if) Suppose that a controller $\mathcal{C}^{\prime}$
satisfies (i) and (ii) above. This controller can be written as
the kernel of a matrix (not necessarily minimal)
$C^{\prime}(\frac{d}{dt})$ with as many rows as
$C(\frac{d}{dt}).$ We know that there is a unimodular matrix $U$ such that%
\begin{equation}
U\left[
\begin{array}
[c]{c}%
R\\
C
\end{array}
\right]  =\left[
\begin{array}
[c]{cc}%
U_{11} & U_{12}\\
U_{21} & U_{22}%
\end{array}
\right]  \left[
\begin{array}
[c]{c}%
R\\
C
\end{array}
\right]  =\left[
\begin{array}
[c]{c}%
R\\
C^{\prime}%
\end{array}
\right]  . \label{6}%
\end{equation}
We shall prove that we can assume $U$ to be of the form%
\begin{equation}
U=\left[
\begin{array}
[c]{cc}%
I & 0\\
V & I
\end{array}
\right]  .
\end{equation}
First, we find a unimodular matrix $W$ such that%
\begin{equation}
RW=\left[
\begin{array}
[c]{cc}%
D & 0
\end{array}
\right]  ,
\end{equation}
where $D$ is a square nonsingular matrix. We then use the following notation%
\begin{align}
\left[
\begin{array}
[c]{c}%
R\\
C
\end{array}
\right]  W  &  =:\left[
\begin{array}
[c]{cc}%
D & 0\\
C_{1} & C_{2}%
\end{array}
\right]  ,\\
\left[
\begin{array}
[c]{c}%
R\\
C^{\prime}%
\end{array}
\right]  W  &  =:\left[
\begin{array}
[c]{cc}%
D & 0\\
C_{1}^{\prime} & C_{2}^{\prime}%
\end{array}
\right]  .
\end{align}
It follows that (\ref{6}) can be rewritten as%
\begin{equation}
U\left[
\begin{array}
[c]{cc}%
D & 0\\
C_{1} & C_{2}%
\end{array}
\right]  W^{-1}=\left[
\begin{array}
[c]{cc}%
D & 0\\
C_{1}^{\prime} & C_{2}^{\prime}%
\end{array}
\right]  W^{-1},
\end{equation}
and since $W$ is unimodular,%
\begin{equation}
U\left[
\begin{array}
[c]{cc}%
D & 0\\
C_{1} & C_{2}%
\end{array}
\right]  =\left[
\begin{array}
[c]{cc}%
U_{11} & U_{12}\\
U_{21} & U_{22}%
\end{array}
\right]  \left[
\begin{array}
[c]{cc}%
D & 0\\
C_{1} & C_{2}%
\end{array}
\right]  =\left[
\begin{array}
[c]{cc}%
D & 0\\
C_{1}^{\prime} & C_{2}^{\prime}%
\end{array}
\right]  .
\end{equation}
Consequently, we have the following equations
\begin{subequations}
\label{7}%
\begin{align}
U_{11}D+U_{12}C_{1}  &  =D,\label{7a}\\
U_{12}C_{2}  &  =0,\label{7b}\\
U_{21}D+U_{22}C_{1}  &  =C_{1}^{\prime},\label{7c}\\
U_{22}C_{2}  &  =C_{2}^{\prime}. \label{7d}%
\end{align}
Since the controller $\mathcal{C}$ is regular, $C_{2}$ must be
full row rank. Now, (\ref{7b}) implies that $U_{12}$ is a left
annihilator of $C_{2}.$
Consequently%
\end{subequations}
\begin{equation}
U_{12}=0.
\end{equation}
Substituting this to (\ref{7a}) yields%
\begin{equation}
U_{11}=I.
\end{equation}
Since $U$ is unimodular, this implies that $U_{22}$ is unimodular.
Thus, we
can conclude that%
\begin{equation}
U=\left[
\begin{array}
[c]{cc}%
I & 0\\
U_{21} & U_{22}%
\end{array}
\right]  ,
\end{equation}
with $U_{22}$ unimodular. Furthermore,
$C^{\prime\prime}:=U_{22}C^{\prime}$ is also a kernel
representation of $\mathcal{C}^{\prime}$ so we can assume $U_{22}$
to be the identity matrix without any loss of generality.
\end{proof}

Now we can parametrize the elements of $\mathfrak{C}_{\mathrm{can}%
}^{\mathrm{reg}}$ as follows.

\begin{theorem}
\label{t4}Let the control manifest behavior of the plant
$\mathcal{P}_{c}$ be the kernel of $P_{c}\left(
\frac{d}{dt}\right)  $ and a controller
$\mathcal{C}\in\mathfrak{C}_{\mathrm{can}}^{\mathrm{reg}}$ be the
kernel of $C\left(  \frac{d}{dt}\right)  $. Assume that both
$P_{c}$ and $C$ are full row rank. A controller
$\mathcal{C}^{\prime}$ is also an element of
$\mathfrak{C}_{\mathrm{can}}^{\mathrm{reg}}$ if and only if it is
the kernel of $V\left(  \frac{d}{dt}\right)  P_{c}\left(
\frac{d}{dt}\right)  +C\left( \frac{d}{dt}\right)  $ for some
polynomial matrix $V(\xi)$.
\end{theorem}

\begin{proof}
The full plant behavior can be represented by%
\begin{equation}
\left[
\begin{array}
[c]{cc}%
\tilde{R}\left(  \frac{d}{dt}\right)  & \tilde{M}\left(
\frac{d}{dt}\right)
\\
0 & P_{c}\left(  \frac{d}{dt}\right)
\end{array}
\right]  \left[
\begin{array}
[c]{c}%
w\\
c
\end{array}
\right]  =0,
\end{equation}
where $\tilde{R}$ is full row rank. It follows that a controller
$\mathcal{C}^{\prime}$ represented as the kernel of
$C^{\prime}\left(
\frac{d}{dt}\right)  $ is regular if and only if%
\begin{equation}
\mathrm{rank}\left[
\begin{array}
[c]{c}%
P_{c}\\
C^{\prime}%
\end{array}
\right]  =\mathrm{rank~}P_{c}+\mathrm{rank~}C^{\prime}.
\end{equation}
This is equivalent to saying that the interconnection of
$\mathcal{P}_{c}$ and $\mathcal{C}^{\prime}$ is regular.
Therefore, we can apply Lemma \ref{l1} (by replacing $\mathcal{K}$
with $\mathcal{C}_{\mathrm{can}}$ and $\mathcal{P}$ with
$\mathcal{P}_{c}$) and obtain the parametrization of all elements
in $\mathfrak{C}_{\mathrm{can}}^{\mathrm{reg}}$.
\end{proof}

\section{Control with minimal interaction}

\subsection{Problem formulation}

Consider the following definition of irrelevant variables.

\begin{definition}
Let a behavior $\mathfrak{B}$ be given by the kernel representation%
\begin{equation}
R_{1}\left(  \frac{d}{dt}\right)  w_{1}+R_{2}\left(
\frac{d}{dt}\right) w_{2}=0.
\end{equation}
If $R_{1}$ is the zero matrix, then the variables in
$\mathbf{w}_{1}$ are said to be \emph{irrelevant} to
$\mathfrak{B}.$
\end{definition}

Notice that whether or not some variables are irrelevant to a
behavior is not a matter of representation. Rather, it is a
property of the behavior. It means for every
$(w_{1},w_{2})\in\mathfrak{B}$ we can always replace $w_{1}$ by
any infinitely differentiable trajectory $w_{1}^{\prime}$ and have
that $(w_{1}^{\prime},w_{2})$ is still an element of
$\mathfrak{B}$. Hence, although $\mathbf{w}_{1}$ is explicitly
present in the description of $\mathfrak{B}$, the information
about its trajectory is irrelevant.

The problem of control with minimal interaction that we are
addressing in this paper can be formulated as follows.\medskip

\noindent\textbf{Control with minimal interaction.} Given are the
full plant behavior $\mathcal{P}_{\mathrm{full}}$ (\ref{1}) and
specification $\mathcal{S}$. We assume that the specification
$\mathcal{S}$ is regularly implementable. Construct a regular
controller $\mathcal{C}$ that implements $\mathcal{S}$ with as
many irrelevant variables as possible, or equivalently as few
relevant variables as possible.\medskip

\begin{figure}
[ptb]
\begin{center}
\includegraphics[
height=3.0709cm, width=8.2703cm
]%
{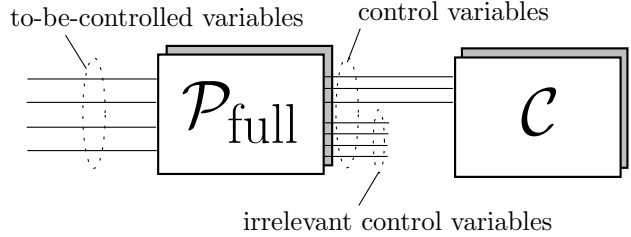}%
\caption{Control with irrelevant control variables.}%
\label{f2}%
\end{center}
\end{figure}

The controller to be designed is called the \emph{controller with
minimal interaction}. When some control variables are irrelevant
to the controller, we can realize the controller without using
these variables. See Figure \ref{f2}
for an illustration.%

\subsection{The solution}

We are going to use the parametrization of $\mathfrak{C}_{\mathrm{can}%
}^{\mathrm{reg}}$ that we derived in the previous section to solve
the problem of control with minimal interaction. First, consider
the following lemma.

\begin{lemma}
\label{l2}Let a behavior $\mathfrak{B}$ be given by the kernel representation%
\begin{equation}
R_{1}\left(  \frac{d}{dt}\right)  w_{1}+R_{2}\left(
\frac{d}{dt}\right) w_{2}=0.
\end{equation}
If $\mathbf{w}_{1}$ is irrelevant to $\mathfrak{B}$, then it is
also irrelevant to any
$\mathfrak{B}^{\prime}\supseteq\mathfrak{B}$.
\end{lemma}

\begin{proof}
The kernel representation of any
$\mathfrak{B}^{\prime}\supseteq\mathfrak{B}$
can be written as%
\[
FR_{1}\left(  \frac{d}{dt}\right)  w_{1}+FR_{2}\left(
\frac{d}{dt}\right) w_{2}=0,
\]
for some polynomial matrix $F$. Clearly $FR_{1}=0$, thus
$\mathbf{w}_{1}$ is irrelevant to $\mathfrak{B}^{\prime}$.
\end{proof}

Lemma \ref{l2} and Theorem \ref{t3} tell us that it is sufficient
to search
for the controller with minimal interaction in $\mathfrak{C}_{\mathrm{can}%
}^{\mathrm{reg}}$, instead of in the set of all regular
controllers. This is
an advantage, since we can parametrize all the controllers in $\mathfrak{C}%
_{\mathrm{can}}^{\mathrm{reg}}$, as shown in Theorem \ref{t4}. To
solve the problem of control with minimal interaction, we need to
find an element of $\mathfrak{C}_{\mathrm{can}}^{\mathrm{reg}}$
with as many zero columns as possible. Generally, since there are
finitely many columns, there is a maximal number of zero columns
that can be attained. However, there is no guarantee that this
number is attained by a unique controller. In fact, generally
speaking, it is not.

The procedure to compute a regular controller that implements
$\mathcal{S}$ and has as many irrelevant variables as possible can
be summarized as follows.

\begin{description}
\item[Step 1.] Construct the canonical controller $\mathcal{C}%
_{\text{\textrm{can}}}$ for the problem. Since $\mathcal{S}$ is
regularly implementable, we know that the canonical controller
implements $\mathcal{S}$.

\item[Step 2.] Construct a controller $\mathcal{C}\in\mathfrak{C}%
_{\mathrm{can}}^{\mathrm{reg}}$. The proof of Theorem \ref{t2}
describes how to construct $\mathcal{C}$ from a regular
controller. Denote the kernel representation of $\mathcal{C}$ and
the control manifest behavior, $\mathcal{P}_{c}$, by
$C(\frac{d}{dt})$ and $P(\frac{d}{dt})$ respectively.

\item[Step 3.] The kernel representation of the controller with
minimal interaction can be found by finding a matrix $V$ such that
$C+VP$ has as many zero columns as possible.
\end{description}

The algebraic problem related to the third step has a
combinatorial aspect in it, as we generally need to search for the
answer by trying all possible subsets of the columns. This
situation gives rise to a computational challenge, namely to
design an algorithm that can handle this combinatorial problem
efficiently. Before we proceed to discuss the algorithm, we
establish an upper bound for the number of irrelevant variables
that can be attained in the controller with minimal interaction.

\begin{lemma}
\label{l4}The controller with minimal interaction can have at most
$\mathtt{c}-\mathtt{p}(\mathcal{C})$ irrelevant variables. Here
$\mathtt{c}$ denotes the number of all control variables (the
cardinality of $\mathbf{c}$) and $\mathtt{p}(\mathcal{C})$ denotes
the number of output variables in $\mathcal{C}$, which is any
regular controller that implements $\mathcal{S}$.
\end{lemma}

\begin{proof}
From the definition of regularity, we know that all regular
controllers that
implement $\mathcal{S}$ have the same number of outputs, i.e., $\mathtt{p}%
(\mathcal{C})$. This is the number of rows in a minimal kernel
representation of the controller. It is easily seen that the
number of columns is
$\mathtt{c}$. If a regular controller has more than $\mathtt{c}-\mathtt{p}%
(\mathcal{C})$ irrelevant variables, then the nonzero entries of
its kernel representation form a tall matrix\footnote{A tall
matrix is a matrix, in which there are more nonzero rows than
there are columns.}, and thus cannot be minimal.
\end{proof}

\begin{notation}
In the subsequent discussion, we denote the entry on the $i$-th
row, $j$-th column of $C$ as $C_{ij}$. The $j$-th column of $C$ is
denoted as $C_{\bullet j}$ and the $i$-th row as $C_{i\bullet}$.
In a similar fashion, we also define $V_{ij}$, $V_{\bullet j}$,
$V_{i\bullet}$, $P_{ij}$, $P_{\bullet j}$, and $P_{i\bullet}$.
Moreover, we denote the greatest common divisor of the polynomials
in $P_{\bullet j}$ as $\pi_{j}$.
\end{notation}

Notice that the $j$-th column of $C+VP$ is zero if and only if
$C_{\bullet j}+VP_{\bullet j}=0$. Consider the following
proposition.

\begin{proposition}
There exists a $V$ such that%
\begin{equation}
C_{\bullet j}+VP_{\bullet j}=0 \label{8}%
\end{equation}
if and only if $\pi_{j}$ divides $C_{\bullet j}$.
\end{proposition}

\begin{proof}
(if) Suppose that $C_{ij}=\kappa_{i}\pi_{j}$, where $\kappa_{i}$
is a polynomial. Since $\pi_{j}$ is the greatest common divisor of
the polynomials in $P_{\bullet j}$, there exists a row vector $v$
such that the Bezout
identity%
\[
v\cdot P_{\bullet j}=\pi_{j}%
\]
is satisfied. It follows that choosing $V$ such that
$V_{ij}=-\kappa_{i}v$ will give us (\ref{8}).

(only if) Suppose that $C_{\bullet j}=-VP_{\bullet j}$. It means
$C_{ij}=-V_{i\bullet}P_{\bullet j}$. Since $\pi_{j}$ divides
$P_{\bullet j}$, it also divides $C_{ij}$.
\end{proof}

The process of constructing the matrix $V$ that corresponds to a
controller with minimal interaction can be considered as a
recursive process.\medskip

\noindent\textbf{The recursive computation for} $V$. Suppose that
we are given $C$ and $P$, and we want to construct $V$ such that
$C+VP$ has as many zero columns as possible. Suppose that first we
want to nullify the $i$-th column of $C+VP$. We proceed with the
following steps:\newline

1. Check if $C_{\bullet i}$ is divisible by $\pi_{i}$. If not, the
procedure stops here, otherwise, we denote $C_{\bullet i}$ as
$\left[
\begin{array}
[c]{cccc}%
\pi_{i}\kappa_{1i} & \pi_{i}\kappa_{2i} & \cdots & \pi_{i}\kappa
_{\mathtt{p}(\mathcal{C})i}%
\end{array}
\right]  ^{T}$.

2. Compute a unimodular matrix $U$ such that $\tilde{P}:=UP$ is
such that its $i$-th column is $\left[
\begin{array}
[c]{cccc}%
\pi_{i} & 0 & \cdots & 0
\end{array}
\right]  ^{T}$.

3. Define $\tilde{V}:=VU^{-1}$. We then have that%
\[
C+VP=C+\tilde{V}\tilde{P}.
\]
It follows that the $i$-th column of $C+VP$ is zero if and only if
the first column of $\tilde{V}$ is $\left[
\begin{array}
[c]{cccc}%
-\kappa_{1i} & -\kappa_{2i} & \cdots & -\kappa_{\mathtt{p}(\mathcal{C})i}%
\end{array}
\right]  ^{T}$.

4. We can write $\tilde{V}$ and $\tilde{P}$ as%
\[
\tilde{V}=\left[
\begin{array}
[c]{cc}%
\tilde{V}_{1} & \tilde{V}_{2}%
\end{array}
\right]  ,\ \tilde{P}=\left[
\begin{array}
[c]{c}%
\tilde{P}_{1}\\
\tilde{P}_{2}%
\end{array}
\right]  ,
\]
where $\tilde{P}_{1}$ is the first row of $\tilde{P}$ and
$\tilde{V}_{1}$ is the first column of $\tilde{V}$, which is now
known.

5. Notice that%
\begin{equation}
C+\tilde{V}\tilde{P}=C+\tilde{V}_{1}\tilde{P}_{1}+\tilde{V}_{2}\tilde{P}_{2}.
\end{equation}
Define $\tilde{C}:=C+\tilde{V}_{1}\tilde{P}_{1}$. The $i$-th
column of $\tilde{C}$ is zero by the construction in the previous
steps. If we want to proceed by, say, nullifying the $j$-th column
of $C+VP$, then the problem is to design $\tilde{V}_{2}$ such that
the $j$-th column of $\tilde{C}+\tilde {V}_{2}\tilde{P}_{2}$ is
zero. So now we arrived at a problem similar to the one we started
with. However, now we have (at least) one less column to nullify
and one less column of $\tilde{V}$ to design (since the first
column is determined). We define an algorithmic function that does
the computation described above.

\noindent\textbf{Algorithm of the function }$(\tilde{C},\tilde{P}_{2}%
,\tilde{V}_{1},U,\mathtt{fail},\mathtt{skip})=\mathtt{nullify}(C,P,i)$

\begin{enumerate}
\item Check if $C_{\bullet i}$ is zero. If yes, return
$(C,P,0,I,\mathtt{false},\mathtt{true})$. Otherwise, go to step 2.

\item Compute a unimodular matrix $U$ such that $\tilde{P}:=UP$ is
such that its $i$-th column is $\left[
\begin{array}
[c]{cccc}%
\pi_{i} & 0 & \cdots & 0
\end{array}
\right]  ^{T}$.

\item Check if $\pi_{i}$ divides $C_{\bullet i}$. If not, $\mathtt{fail}%
=\mathtt{true}$, return $(0,0,0,U,\mathtt{fail})$. Otherwise $\mathtt{fail}%
=\mathtt{false}$ and go to step 4.

\item Compute $\tilde{V}_{1}=-\frac{1}{\pi_{i}}C_{\bullet i}$.

\item Partition $\tilde{P}$ into $\left[
\begin{array}
[c]{c}%
\tilde{P}_{1}\\
\tilde{P}_{2}%
\end{array}
\right]  ,$ where $\tilde{P}_{1}$ is the first row of $\tilde{P}$.

\item Compute $\tilde{C}=C+\tilde{V}_{1}\tilde{P}_{1}$.

\item Return $(\tilde{C},\tilde{P}_{2},\tilde{V}_{1},U,\mathtt{fail}%
,\mathtt{false})$.
\end{enumerate}

The following example illustrates the algorithm.

\begin{example}
\label{ex3}Let%
\[
P(\xi)=\left[
\begin{array}
[c]{ccc}%
\xi & \xi & 1\\
\xi+1 & \xi & 0
\end{array}
\right]  ,C(\xi)=\left[
\begin{array}
[c]{ccc}%
\xi^{2}-\xi & \xi-1 & -1
\end{array}
\right]  .
\]
Suppose that we want to design $V(\xi)$ such that the first column
of $C+VP$ is zero, i.e. $i=1$. We start with step 1 of the
algorithm above, and since the first column of $C$ is not zero,
\texttt{skip=false} and we go to step 2,
where we obtain%
\begin{align*}
U(\xi) &  =\left[
\begin{array}
[c]{cc}%
-1 & 1\\
-\xi-1 & \xi
\end{array}
\right]  ,\\
\tilde{P}(\xi) &  :=U(\xi)P(\xi)=\left[
\begin{array}
[c]{ccc}%
1 & 0 & -1\\
0 & -\xi & -\xi-1
\end{array}
\right]  .
\end{align*}
Notice that the greatest common divisor of $P_{\bullet i}$ is 1,
which divides $C_{\bullet i}.$ This means
$\mathtt{fail}=\mathtt{false}$ and we go to step
4, 5 and 6.%
\begin{align*}
\tilde{V}_{1}(\xi) &  =\xi-\xi^{2},~\tilde{P}_{1}(\xi)=\left[
\begin{array}
[c]{ccc}%
1 & 0 & -1
\end{array}
\right]  ,\\
\tilde{P}_{2}(\xi) &  =\left[
\begin{array}
[c]{ccc}%
0 & -\xi & -\xi-1
\end{array}
\right]  ,~\tilde{C}(\xi)=\left[
\begin{array}
[c]{ccc}%
0 & \xi-1 & \xi^{2}-\xi-1
\end{array}
\right]
\end{align*}
Thus, we can verify that the first column of $\tilde{C}$ has been
nullified.
\end{example}

\begin{figure}
[ptb]
\begin{center}
\includegraphics[
height=0.9997in, width=1.3111in
]%
{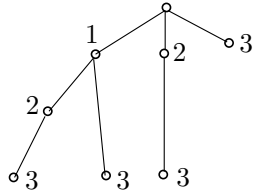}%
\caption{The search tree for $\mathtt{c}=3$.}%
\label{f3}%
\end{center}
\end{figure}

The decision on the order of the columns that we nullify involves
a combinatorial search. Now, we are going to develop a depth-first
search algorithm\footnote{Depth-first search is a standard term in
computer science. It is an algorithm for traversing or searching a
tree, tree structure, or graph. Intuitively, in this algorithm,
one starts at the root (selecting some node as the root in the
graph case) and explores as far as possible along each branch
before backtracking.} (see for example, \cite{Cormen}) that can
compute the controller with minimal interaction. First, we are
going to define an efficient data representation for the
depth-first search. Consider the set
$\mathbf{A}:=\{1,2,3,\cdots,\mathtt{c}\}$, where $\mathtt{c}$ is
the cardinality of $\mathbf{c}$, which is also the number of
columns in the representation of the controller. We define
$\mathbf{S}$ to be the set of increasing strings formed by the
elements of $\mathbf{A}$ nonrepeatingly, such that if
$s\in\mathbf{S}$ and $\left\vert s\right\vert $ is the length of
$s$
then%
\begin{equation}
s_{i}<s_{i+1},\ 1\leq i<\left\vert s\right\vert .
\end{equation}
The empty string is denoted by $\varepsilon$ and its length is
defined to be 0. Thus, the numbers in a string $s\in\mathbf{S}$ in
increasing. There are $2^{\mathtt{c}}$ elements of $\mathbf{S}$,
each of which represents an element of the power set of
$\mathbf{A}$. We can assign lexicographic ordering (see for
example, \cite{Cormen}) on the elements of $\mathbf{S}$ and sort
them.
That is, for any two distinct elements $s,s^{\prime}\in\mathbf{S}$,%
\begin{equation}
s<s^{\prime}:\Leftrightarrow\sum_{i=1}^{\mathtt{c}}(s_{i}^{\prime}-s_{i}%
)N^{i}>0,
\end{equation}
where $N$ is any integer larger than $\mathtt{c}$, and under the
convention
that%
\begin{equation}
s_{i}=0,i>\left\vert s\right\vert .
\end{equation}
We then define a subset $\mathbf{T}\subset\mathbf{S}$, by%
\begin{equation}
\mathbf{T}:=\{s\in\mathbf{S}~|~s_{\left\vert s\right\vert
}=\mathtt{c}\}.
\end{equation}
That is, $\mathbf{T}$ is the set of strings that end with
$\mathtt{c}$. For example, the elements of $\mathbf{S}$ and
$\mathbf{T}$ for $\mathtt{c}=3$, in ascending order, are
\begin{align*}
\mathbf{S}  &  =\left\{  \varepsilon,1,12,123,13,2,23,3\right\}  ,\\
\mathbf{T}  &  =\{123,13,23,3\}.
\end{align*}
The search tree for the problem, where $\mathtt{c}=3$ is shown in
Figure \ref{f3}. Notice that each element of $\mathbf{S}$
represents a node on this tree. The ordering of $\mathbf{S}$ tells
us the order in which the nodes are visited in the depth-first
search. The elements of $\mathbf{T}$ represent the terminal nodes,
each which represents a path from the initial node (the top of
the graph) to a terminal node.%

We define the following operations on $\mathbf{S}$. The prefix
operator $pre:\mathbf{S}\rightarrow2^{\mathbf{S}}$ is such that
$pre(s)$ is the set containing all the prefixes of $s$. The
operator $\bullet^{+}$ and $\bullet^{-}$ are such that $s^{+}$ is
the last symbol in the string $s$ and $s^{-}$ is the string formed
by removing the last symbol from $s$. The
operator $\left\vert \cdot,\cdot\right\rangle :\mathbf{S}\times\mathbb{N}%
\rightarrow\mathbf{T}$ is defined as follows.%
\begin{equation}
\left\vert s,k\right\rangle :=\text{the smallest element of
}\{\sigma
\in\mathbf{T~|~}\sigma>s,\left\vert \sigma\right\vert >k\}. \label{9}%
\end{equation}
If the set in (\ref{9}) is empty, then $\left\vert
s,k\right\rangle :=\varepsilon.$ The operation $\left\vert
\cdot,\cdot\right\downarrow
:\mathbf{S}\times\mathbf{S}\rightarrow\mathbf{S}$ is defined as follows.%
\begin{equation}
\left\vert s,s^{\prime}\right\downarrow :=\text{the smallest
element of }pre(s^{\prime})\cap pre(s)^{\complement}.
\end{equation}
The operation $\left\lceil \cdot\right\rceil
:\mathbf{S}\rightarrow\mathbf{S}$
is defined as follows.%
\begin{equation}
\left\lceil s\right\rceil :=\text{the smallest
}s^{\prime}\in\mathbf{S}\text{ s.t. }s\notin pre(s^{\prime}).
\end{equation}
Denote the cardinality of $\mathbf{S}$ by $\left\vert
\mathbf{S}\right\vert =2^{\mathtt{c}}-1$, and the $i$-th element
of $\mathbf{S}$ by $\mathbf{S}(i)$.

The following algorithm takes polynomial matrices $C$ and $P$ as
inputs and returns a matrix $V$ such that $C+VP$ has as many zero
columns as possible.

\noindent\textbf{Algorithm of the function
}$V=\mathtt{computeV}(C,P)$

\begin{enumerate}
\item If $C$ already has zero columns, compute a unimodular
permutation $W_{1}$ such that $\hat{C}:=CW_{1}$ has all the zero
columns on the left. Otherwise, $W_{1}=I$. Define
$\hat{P}:=PW_{1}$.

\item Denote the maximum possible number of zero columns (see
Lemma \ref{l4}) as $M.$ If the number of zero columns of $\hat{C}$
is less than $M$, go to step 3, otherwise return $V=0$.

\item Initialize the variable $k=0,$ $n_{\max}=0,$ and the strings
$\sigma$ and $\sigma_{\max}$ are both empty.

\item Define $\tilde{C}[0]=\hat{C}$, $\tilde{P}[0]=\hat{P},$ $\mathtt{fail}%
(0)=\mathtt{false}$.

\item If $n_{\max}<M$ and $k<\left\vert \mathbf{S}\right\vert $,
go to step 6. Otherwise, go to step 13.

\item Define $\sigma=\mathbf{S}(k).$

\item If $\mathtt{fail}(k)=\mathtt{true}$ then change $k$ such
that $\mathbf{S}(k)=\left\lceil \sigma\right\rceil $ and go to
step 9, otherwise go to step 8.

\item If $\{s\in\mathbf{T~|~}s>\sigma,\left\vert s\right\vert
>n_{\max
}\}=\emptyset$ go to step 13, otherwise change $k$ such that $\mathbf{S}%
(k)=\left\vert \sigma,\left\vert \sigma,n_{\max}\right\rangle
\right\downarrow $.

\item Let $j$ be such that $\mathbf{S}(j)=\left(
\mathbf{S}(k)\right)  ^{-}$, compute
$(\tilde{C}[k],\tilde{P}_{2}[k],$ $\tilde{V}_{1}[k],U[k],$
$\mathtt{fail}[k],\mathtt{skip}[k])=\mathtt{nullify}(C[j],P[j],(\mathbf{S}%
(k))^{+})$.

\item If $\mathtt{fail}[k]=\mathtt{false}$, then go to step 11,
otherwise go to step 12.

\item If $\left\vert \mathbf{S}(k)\right\vert >n_{\max}$ then
modify $n_{\max }=\left\vert \mathbf{S}(k)\right\vert $ and
$\sigma_{\max}=\mathbf{S}(k)$.

\item Go to step 5.

\item Initialize $\sigma=\sigma_{\max}$. If
$\sigma_{\max}=\varepsilon$, then return $V=0$. Otherwise go to
step 14.

\item If $\left\vert \sigma\right\vert =0$ go to step 19,
otherwise go to step 15.

\item Let $k$ be such that $\mathbf{S}(k)=\sigma$. If $\mathtt{skip}%
[k]=\mathtt{true}$ then go to step 18. Otherwise go to step 16.

\item
$V=0_{\operatorname{row}\dim(C)\times\operatorname{col}\dim(\tilde
{P}_{2}[k])}$.

\item $\tilde{V}=\left[
\begin{array}
[c]{cc}%
\tilde{V}_{1}[k] & V
\end{array}
\right]  ,$ redefine $V=\tilde{V}U[k].$

\item $\sigma=\left(  \mathbf{S}(k)\right)  ^{-}$ and go to step
14.

\item Return $V$.
\end{enumerate}

\begin{example}
Consider the matrices given in Example \ref{ex3}.%
\[
P(\xi)=\left[
\begin{array}
[c]{ccc}%
\xi & \xi & 1\\
\xi+1 & \xi & 0
\end{array}
\right]  ,C(\xi)=\left[
\begin{array}
[c]{ccc}%
\xi^{2}-\xi & \xi-1 & -1
\end{array}
\right]  .
\]
If we apply the algorithm above to this example, then the
following steps are going to executed.\newline1. The first column
of $C$ will be nullified, as shown in Example \ref{ex3}. Thus the
maximum number of column that can be nullified by the algorithm so
far is 1.\newline2. See Figure \ref{f3}. The algorithm is now at
the first branch from the left, in the figure. The algorithm will
subsequently try to nullify the second column. That is, it will
try to find polynomial matrix $\tilde{V}_{2}$ such that the second column of%
\[
\left[
\begin{array}
[c]{ccc}%
0 & \xi-1 & \xi^{2}-\xi-1
\end{array}
\right]  +\tilde{V}_{2}\left[
\begin{array}
[c]{ccc}%
0 & -\xi & -\xi-1
\end{array}
\right]
\]
is zero. Since this is not possible, it will try to nullify the
third column, which corresponds to the second branch from the left
in Figure \ref{f3}. This is also not possible.\newline3. The
algorithm then tries to nullify the second column of $C$ (third
branch from the left in Figure \ref{f3}). This is not possible
since the greatest common divisor of $P_{\bullet2}$ is $\xi$, and
it does not divide $C_{\bullet2}$.\newline4. Since the remaining
branch in Figure \ref{f3} consists of only one element, it is not
possible to find a combination of columns, consisting of more than
one column, that can be nullified. The algorithm then terminates
and the final result is the nullification of the first column of
$C$.\newline
\end{example}

\section{Control problem with input-output partition constraint}

\subsection{Problem formulation}

One of the features of the behavioral approach to systems theory
is that no \emph{a priori} distinction is made between input and
output variables of a system \cite{Willems91,Willems98}. This
means that given a certain law that describes the system, the
system is identified by the collection of its trajectories as is.
Therefore, it is not necessary to have any input-output structure
when describing the system.

However, when two systems are interconnected, sometimes some
input-output
structure can emerge naturally as a constraint. Consider the following example.%

\begin{figure}
[ptb]
\begin{center}
\includegraphics[
height=7.7958cm, width=8.1231cm
]%
{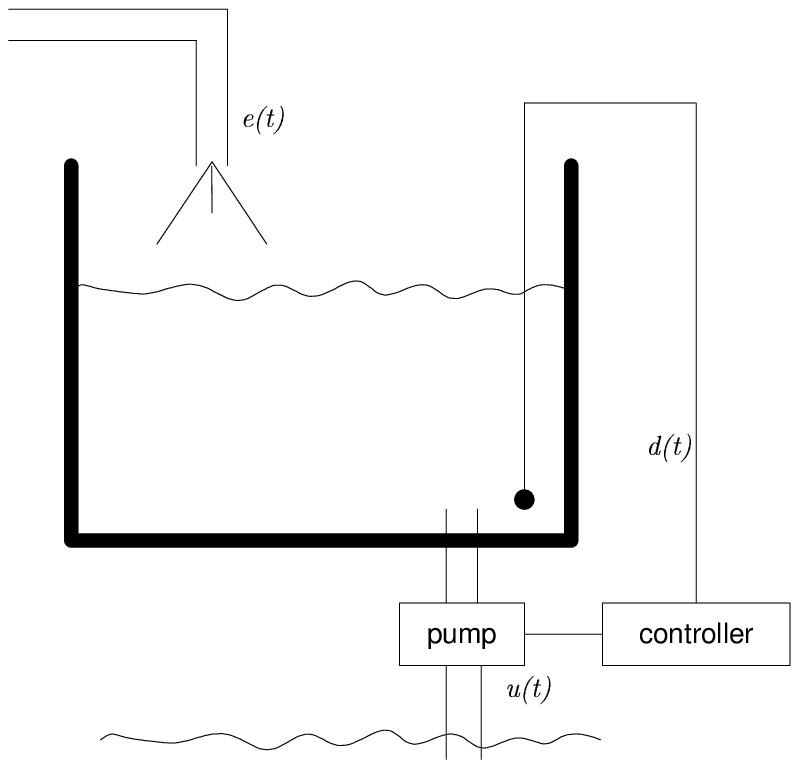}%
\caption{The water tank system in Example \ref{ex1}.}%
\label{f4}%
\end{center}
\end{figure}

\begin{example}
\label{ex1}Consider a tank filled with water as shown in Figure
\ref{f4}. On top of the tank is an inlet from which a variable
flow of water can get into the tank. We denote the water flow from
this inlet as $\mathbf{e}.$ On the bottom of the tank, there is an
opening connected to a pump that can pump water out of/into the
tank. We denote the amount of water flow pumped out of the tank as
$\mathbf{u}$. The tank is also equipped with a sensor that
measures the change of volume of water inside the tank, the
measurement of the sensor is denoted as $\mathbf{d}$. The
mathematical model of this system can
be simply written as%
\begin{equation}
d(t)=e(t)-u(t).
\end{equation}
Now consider the following control problem. Given $\mathbf{d}$ and
$\mathbf{u}$ as control variables, we want to design a controller
such that the level of water is constant, i.e. $e(t)=u(t).$ In
other words, we aim at perfect tracking of $\mathbf{e}$ by
$\mathbf{u}$. Intuitively, we know that such task cannot be
accomplished. However, consider the following
construction. First we write the plant behavior in a kernel representation.%
\begin{equation}
\mathcal{P}=\{(e,u,d)~|~e(t)-u(t)-d(t)=0\}. \label{10}%
\end{equation}
We then take a candidate controller $\mathcal{C}$ expressed by%
\begin{equation}
\mathcal{C}=\{(u,d)~|~d(t)=0\}. \label{11}%
\end{equation}
The interconnection $\mathcal{P}\parallel\mathcal{C}$ is
represented by the
\begin{equation}
\left[
\begin{array}
[c]{ccc}%
1 & -1 & -1\\
0 & 0 & 1
\end{array}
\right]  \left[
\begin{array}
[c]{c}%
e\\
u\\
d
\end{array}
\right]  =0. \label{12}%
\end{equation}
Notice that the interconnection exhibits the following
features.\newline(i) The interconnection is a regular
interconnection. In fact, it is even a regular feedback
interconnection \cite{Willems97}.\newline(ii) The controller is
expressed only in terms of $\mathbf{u}$ and
$\mathbf{d}$.\newline(iii) In the controlled behavior, perfect
tracking $e(t)=u(t)$ is attained.\newline
\end{example}

In the example above, the proposed controller is regular and
accomplishes the control task. However, this is still counter
intuitive, and impossible to implement. The variable $\mathbf{d}$
is a measurement coming from a sensor, and yet we use it to
enforce control on the system. Otherwise stated, we control the
system by restricting the reading of a sensor. Now, consider the
following modification of the example.

\begin{example}
\label{ex2}Let us swap the name of variables involved in the
system as follows. We swap $\mathbf{d}$ and $\mathbf{u}$. The
schematic of the system is
now shown in Figure \ref{f5}.%
\begin{figure}
[ptb]
\begin{center}
\includegraphics[
height=3.0692in, width=3.1981in
]%
{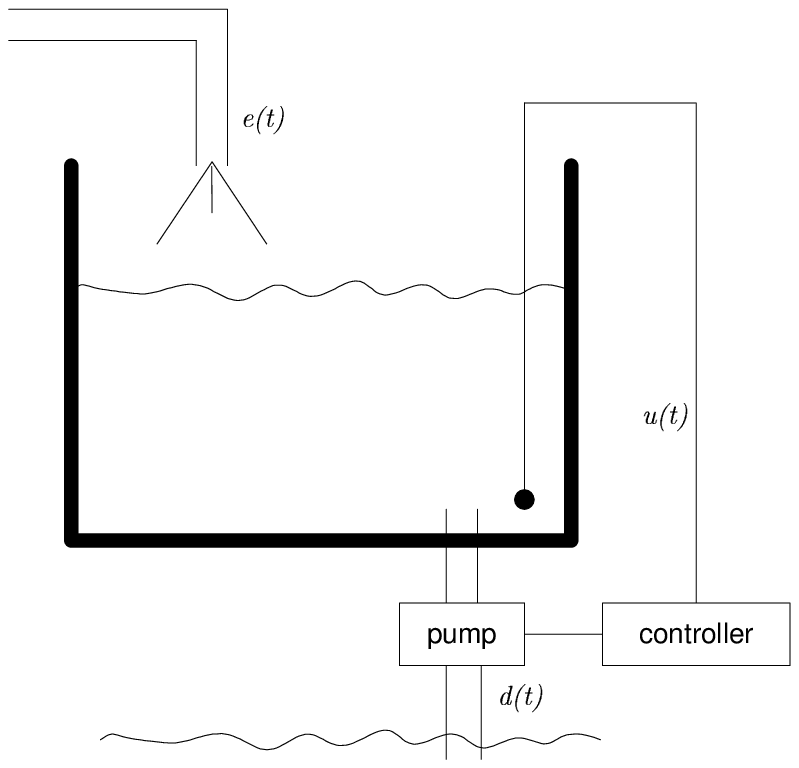}%
\caption{The water tank system in Example \ref{ex2}.}%
\label{f5}%
\end{center}
\end{figure}
Notice that the mathematical model of the system is still given by
(\ref{10}). Now take the controller $\mathcal{C}$ given by
(\ref{11}). Clearly, the features of the interconnection
(\ref{12}) are still there. What the controller now does is shut
down the pump. This controller does not keep the water level
constant. But, that is not the fact that we are interested in. The
interesting observation is that now the interconnection does make
sense.
\end{example}

These two examples suggest the following facts.

\begin{itemize}
\item We may need to introduce a \emph{constraint} for systems
interconnection to make sense. The constraint cannot be formulated
based on the mathematical representation of the systems alone. The
interconnections described in Example \ref{ex1} and Example
\ref{ex2} share the same mathematical representation, yet in one
case the constraint is not satisfied, while in the other it is.
This is in contrast with the regularity constraint, where the
constraint can actually be derived from the behaviors themselves.

\item The new constraint is different from the regularity
constraint. Example \ref{ex1} describes an interconnection where
the regularity constraint is satisfied, while the new constraint
that we are going to formulate is not satisfied.
\end{itemize}

As indicated by the Example \ref{ex1}, the constraint is violated
when the plant is restricted through a variable that is inherently
an output of the system. That is, the variable is physically
dictated to be an output of the system. The information that a
variable is an output cannot be deduced from the mathematical
description of the system, rather it has to be provided in
addition to the description of the plant. System variables that
have to be output variables by physical consideration, are called
\emph{declared outputs}. We then require that\emph{ the controller
accepts the declared output of the plant as its input}, for the
interconnection to make sense. To say it differently, suppose that
$\mathbf{y}$ is a (set of) variable(s) that is a part of the
control variables. If $\mathbf{y}$ is declared as output because
of some physical interpretation of the system, we want to
input-output partition the variables of the controller, such that
$\mathbf{y}$ belongs to the input part. Input-output partitioning
of the variables of a linear system has been introduced in Section
\ref{sec2}.

The control problem with input-output partitioning constraint for
linear systems is then formally defined as follows.

\noindent\textbf{Control with input-output partition constraint.}
Given a
control problem, where the plant is%
\begin{equation}
\mathcal{P}=\left\{  (w,u,y)~|~R\left(  \frac{d}{dt}\right)
w+P\left( \frac{d}{dt}\right)  u+Q\left(  \frac{d}{dt}\right)
y=0\right\}  .
\end{equation}
The control variables are $\mathbf{u}$ and $\mathbf{y}$, where
$\mathbf{y}$ is the declared output variables of the plant. The
to-be-controlled variable is
$\mathbf{w}$. The desired specification is given as%
\begin{equation}
\mathcal{S}=\left\{  w~|~S\left(  \frac{d}{dt}\right)  w=0\right\}
.
\end{equation}
Find a regular controller $\mathcal{C}$ described as%
\begin{equation}
\mathcal{C}=\left\{  (u,y)~|~C_{1}\left(  \frac{d}{dt}\right)
u+C_{2}\left( \frac{d}{dt}\right)  y=0\right\}  ,
\end{equation}
such that $\mathcal{C}$ implements $\mathcal{S}$ and the variables
in $\mathcal{C}$ can be input-output partitioned such that
$\mathbf{y}$ belongs to the input part.

\subsection{The solution}

We shall now devise an algorithm that solves the problem. We
assume that the specification $\mathcal{S}$ is regularly
achievable (otherwise the problem is clearly not solvable)

\begin{notation}
We denote the class of regular controllers that implements
$\mathcal{S}$ as $\mathfrak{C}_{\mathcal{S}}^{\mathrm{{reg}}}$.
\end{notation}

To find a solution to the problem, we need to use the following
result.

\begin{lemma}
Given a controller%
\begin{equation}
\mathcal{C}=\left\{  (u,y)~|~C_{1}\left(  \frac{d}{dt}\right)
u+C_{2}\left( \frac{d}{dt}\right)  y=0\right\}  .
\end{equation}
Without loss of generality we assume that $[C_{1}\ C_{2}]$ is full
row rank. The following statements are equivalent.\newline(i) The
variables in $\mathcal{C}$ can be partitioned such that
$\mathbf{y}$ belongs to the input part\newline(ii) $C_{1}$ is full
row rank.\newline(iii) For any $y\in
\mathfrak{C}^{\infty}(\mathbb{R},\mathbb{R}^{\mathtt{y}})$ there
exists a
$u\in\mathfrak{C}^{\infty}(\mathbb{R},\mathbb{R}^{\mathtt{u}})$
such that $(u,y)\in\mathcal{C}$.
\end{lemma}

\begin{proof}
(ii $\Rightarrow$ i) Suppose that $C_{1}$ is full row rank. If
$C_{1}$ is a square matrix, then we already have an input-output
partition with $\mathbf{u}$ as the output and $\mathbf{y}$ as the
input. If $C_{1}$ is not
square, then we can partition it into%
\begin{equation}
C_{1}=\left[
\begin{array}
[c]{cc}%
C_{11} & C_{12}%
\end{array}
\right]  , \label{3.27}%
\end{equation}
possibly after rearranging the columns, such that $C_{11}$ is a
square matrix with full row rank. We can also partition
$\mathbf{u}$ accordingly into $\mathbf{u}_{1}$ and
$\mathbf{u}_{2}$. Now we have an input-output partition with
$\mathbf{u}_{1}$ as the output and $\mathbf{u}_{2}$ and
$\mathbf{y}$ as the input.

(i $\Rightarrow$\ iii) Suppose that the variables in $\mathcal{C}$
can be partitioned such that $\mathbf{y}$ belongs to the input
partition. This means we can partition $\mathbf{u}$ into
$\mathbf{u}_{1}$ and $\mathbf{u}_{2},$ such that we have
$\mathbf{u}_{1}$ as the output and $\mathbf{u}_{2}$ and
$\mathbf{y}$ as the input. So we can partition $C$ accordingly
such that (\ref{3.27}) holds. Following the elimination procedure
in Section \ref{sec2}, we can eliminate $\mathbf{u}_{1}$ and find
that the behavior in terms of
$\mathbf{y}$ and $\mathbf{u}_{2}$ is $\mathfrak{C}^{\infty}(\mathbb{R}%
,\mathbb{R}^{\mathtt{y}+\mathtt{u}_{2}}).$

(iii $\Rightarrow$\ ii) We shall prove it by contradiction.
Suppose that $C_{1}$ is not full row rank. The matrix $[C_{1}\
C_{2}]$ can be transformed
(by premultiplication with a suitable unimodular matrix) into%
\[
\left[
\begin{array}
[c]{cc}%
C_{1}^{\prime} & C_{21}^{\prime}\\
0 & C_{22}^{\prime}%
\end{array}
\right]  ,
\]
where $C_{1}^{\prime}$ and $C_{22}^{\prime}$ are full row rank.
Following the elimination procedure in Section \ref{sec2}, we can
eliminate $\mathbf{u}$ and find that the behavior in terms of
$\mathbf{y}$ is the kernel of $C_{22}^{\prime}(\frac{d}{dt})$.
Hence, we cannot choose any $y\in
\mathfrak{C}^{\infty}(\mathbb{R},\mathbb{R}^{\mathtt{y}})$ as a
trajectory of $\mathbf{y}$.
\end{proof}

It is straightforward to see that problem can be reformulated as
follows.

\noindent\textbf{Problem. }Find a controller $\mathcal{C}\in\mathfrak{C}%
_{\mathcal{S}}^{\mathrm{{reg}}}$ in the form of%
\begin{equation}
\mathcal{C}=\left\{  (u,y)~|~C_{1}\left(  \frac{d}{dt}\right)
u+C_{2}\left(
\frac{d}{dt}\right)  y=0\right\}  , \label{14}%
\end{equation}
where $C_{1}$ is full row rank.

We shall use the following result.

\begin{lemma}
Let $X$ be a subset of
$\mathfrak{C}_{\mathcal{S}}^{\mathrm{{reg}}}$ such that for any
$\mathcal{C\in}\mathfrak{C}_{\mathcal{S}}^{\mathrm{{reg}}}$ there
exists a $\mathcal{C}^{\prime}\in X$ such that
$\mathcal{C}\subseteq
\mathcal{C}^{\prime}$. Then there exists a $\mathcal{C}\in\mathfrak{C}%
_{\mathcal{S}}^{\mathrm{{reg}}}$ that solves the control problem
with input-output partitioning constraint if and only if there
exists a $\mathcal{C}^{\prime}\in X$ that does so.
\end{lemma}

\begin{proof}
(if) Trivial, since
$X\subset\mathfrak{C}_{\mathcal{S}}^{\mathrm{{reg}}}$.

(only if) Suppose that $\mathcal{C}\in\mathfrak{C}_{\mathcal{S}}%
^{\mathrm{{reg}}}$ satisfies the constraint. We shall show that
any
$\mathcal{C}^{\prime}\in\mathfrak{C}_{\mathcal{S}}^{\mathrm{{reg}}}$
such that $\mathcal{C}\subseteq\mathcal{C}^{\prime}$ also
satisfies the constraint. Let $\mathcal{C}$ be given as the kernel
of $\left[
\begin{array}
[c]{cc}%
C_{1} & C_{2}%
\end{array}
\right]  $ as in (\ref{14}). We know that $C_{1}$ is full row
rank. Since $\mathcal{C}\subseteq\mathcal{C}^{\prime}$, there must
be a full row rank matrix $F$ such that $\mathcal{C}^{\prime}$ is
the kernel of $\left[
\begin{array}
[c]{cc}%
FC_{1} & FC_{2}%
\end{array}
\right]  $. We also know that $FC_{1}$ is full row rank. Therefore
$\mathcal{C}^{\prime}$ also solves the problem.
\end{proof}

This lemma tells us that if we can construct a subset of $\mathfrak{C}%
_{\mathcal{S}}^{\mathrm{{reg}}}$ with the property of $X$, we do
not need to
search for the candidate controller in the whole $\mathfrak{C}_{\mathcal{S}%
}^{\mathrm{{reg}}}$. Rather, we can restrict our attention in $X$.
Theorem \ref{t3} shows that
$\mathfrak{C}_{\mathrm{can}}^{\mathrm{reg}}$ has the desired
property. Thus, we shall try to construct the desired controller
in $\mathfrak{C}_{\mathrm{can}}^{\mathrm{reg}}$, which we can
parametrize according to Theorem \ref{t4}.

A solution to the control problem can be found by executing the
following steps.

\begin{description}
\item[Step 1.] Construct the canonical controller $\mathcal{C}%
_{\text{\textrm{can}}}$ for the problem. Since $\mathcal{S}$ is
regularly implementable, we know that the canonical controller
implements $\mathcal{S}$.

\item[Step 2.] Construct a controller $\mathcal{C}\in\mathfrak{C}%
_{\mathrm{can}}^{\mathrm{reg}}$. The proof of Theorem \ref{t2}
contains information on how to construct $\mathcal{C}$ from a
regular controller. Denote the kernel representation of
$\mathcal{C}$ and the control manifest
behavior, $\mathcal{P}_{c}$, by%
\begin{align}
\mathcal{C}  &  =\left\{  (u,y)~|~C_{1}\left(  \frac{d}{dt}\right)
u+C_{2}\left(  \frac{d}{dt}\right)  y=0\right\}  ,\\
\mathcal{P}_{c}  &  =\left\{  (u,y)~|~P_{1}\left(
\frac{d}{dt}\right) u+P_{2}\left(  \frac{d}{dt}\right)
y=0\right\}  ,
\end{align}
respectively.

\item[Step 3.] Following Theorem \ref{t4}, any controller $\mathcal{C}%
^{\prime}$ in $\mathfrak{C}_{\mathrm{can}}^{\mathrm{reg}}$ can be
represented
as%
\[
\mathcal{C}^{\prime}=\left\{  (u,y)~|~\left(  C_{1}+VP_{1}\right)
\left(
\frac{d}{dt}\right)  u+\left(  C_{2}+VP_{2}\right)  \left(  \frac{d}%
{dt}\right)  y=0\right\}
\]
The kernel representation of a controller in $\mathfrak{C}_{\mathrm{can}%
}^{\mathrm{reg}}$ that satisfies the input-output partitioning
constraint can be found by finding a matrix $V$ such that
$C_{1}+VP_{1}$ is full row rank.
\end{description}

The necessary and sufficient condition for the existence of such a
matrix $V$ is given in the following lemma.

\begin{lemma}
\label{l3}Given polynomial matrices
$C\in\mathbb{R}^{\mathtt{c}\times \mathtt{q}}[\xi]$ and
$P\in\mathbb{R}^{\mathtt{p}\times\mathtt{q}}[\xi]$.
There exists a polynomial matrix $V\in\mathbb{R}^{\mathtt{c}\times\mathtt{p}%
}[\xi]$ such that $C+VP$ is full row rank if and only if%
\begin{equation}
\text{rank}\left[
\begin{array}
[c]{c}%
P\\
C
\end{array}
\right]  \geq\mathtt{c.} \label{16}%
\end{equation}

\end{lemma}

\begin{proof}
(only if) Consider the following relation%
\begin{equation}
\left[
\begin{array}
[c]{cc}%
I & 0\\
V & I
\end{array}
\right]  \left[
\begin{array}
[c]{c}%
P\\
C
\end{array}
\right]  =\left[
\begin{array}
[c]{c}%
P\\
C+VP
\end{array}
\right]  . \label{15}%
\end{equation}
Suppose that $C+VP$ is full row rank. This means it has a rank of $\mathtt{c}%
$. Since the left multiplication in (\ref{15}) is a unimodular
transformation,
we know that%
\begin{equation}
\text{rank}\left[
\begin{array}
[c]{c}%
P\\
C
\end{array}
\right]  =\text{rank}\left[
\begin{array}
[c]{c}%
P\\
C+VP
\end{array}
\right]  \geq\mathtt{c.}%
\end{equation}

(if) Assume that (\ref{16}) holds. If $C$ or $P$ is zero, we can
obviously choose a $V$ such that $C+VP$ full row rank. We exclude
these trivial cases and suppose that both $P$ and $C$ are nonzero.
Since the rank of a polynomial matrix is not affected by left and
right multiplication by unimodular matrices, we can assume without
any loss of generality that $P$ has the the
Smith form, that is, the form of%
\begin{equation}
P=\left[
\begin{array}
[c]{cc}%
P_{1} & 0\\
0 & 0
\end{array}
\right]  ,
\end{equation}
where $P_{1}$ is a diagonal matrix with nonzero determinant.
Furthermore, with some appropriate left multiplication with
unimodular matrix, we can transform
$C$ to the following form.%
\begin{equation}
C=\left[
\begin{array}
[c]{cc}%
C_{11} & C_{12}\\
C_{21} & 0\\
0 & 0
\end{array}
\right]  ,
\end{equation}
where $C_{12}$ and $C_{21}$ are full row rank. Denote the rank of
$P_{1}$, $C_{12}$, and $C_{21}$ as $\mathtt{p}^{\prime}$,
$\mathtt{c}^{\prime}$, and
$\mathtt{c}^{\prime\prime}$ respectively. We have the following relation%
\begin{align}
\text{rank}\left[
\begin{array}
[c]{c}%
P\\
C
\end{array}
\right]   &  =\text{rank }C_{12}+\text{rank }\left[
\begin{array}
[c]{c}%
P_{1}\\
C_{21}%
\end{array}
\right]  ,\\
&  =\mathtt{c}^{\prime}+\mathtt{p}^{\prime},\\
&  \geq\mathtt{c.}%
\end{align}
Thus%
\begin{equation}
\mathtt{p}^{\prime}\geq\mathtt{c}-\mathtt{c}^{\prime}. \label{17}%
\end{equation}
We can partition $V$ accordingly to form%
\begin{equation}
V=\left[
\begin{array}
[c]{cc}%
V_{11} & V_{12}\\
V_{21} & V_{22}\\
V_{31} & V_{32}%
\end{array}
\right]  .
\end{equation}
We structure $V$ to have the following form.%
\begin{equation}
V=\left[
\begin{array}
[c]{cc}%
0 & 0\\
0 & 0\\
V_{31} & 0
\end{array}
\right]  ,
\end{equation}
where $V_{31}$ is to be chosen later. Therefore%
\begin{equation}
C+VP=\left[
\begin{array}
[c]{cc}%
C_{11} & C_{12}\\
C_{21} & 0\\
V_{31}P_{1} & 0
\end{array}
\right]  .
\end{equation}
Our goal is to make $C+VP$ a full row rank matrix. Therefore,
$\left[
\begin{array}
[c]{c}%
C_{21}\\
V_{31}P_{1}%
\end{array}
\right]  $ has to be a full row rank matrix. Since $C_{21}$ is
full row rank
and has the rank of $\mathtt{c}^{\prime\prime}$, we can find $\mathtt{c}%
^{\prime\prime}$ columns of $C_{21}$ that form a square matrix
with nonzero determinant. Denote this selection as $N$, where
$N\subset\{1,2,\cdots
,\mathtt{p}^{\prime}\}$. We construct $V_{31}\in\mathbb{R}^{(\mathtt{c}%
-\mathtt{c}^{\prime}-\mathtt{c}^{\prime\prime})\times\mathtt{m}^{\prime}}%
[\xi]$ such that the entries on the $i-$th column of $V_{31}$ are
zero if $i\in N$. The remaining
$(\mathtt{m}^{\prime}-\mathtt{c}^{\prime\prime})$
columns of $V_{31}$ form a $(\mathtt{c}-\mathtt{c}^{\prime}-\mathtt{c}%
^{\prime\prime})$ by
$(\mathtt{m}^{\prime}-\mathtt{c}^{\prime\prime})$ matrix. From
(\ref{17}) we know that it is a wide matrix. We choose the values
of the entries of these columns such that this wide matrix is full
row rank. It follows that $\left[
\begin{array}
[c]{c}%
C_{21}\\
V_{31}P_{1}%
\end{array}
\right]  $ is a full row rank matrix and hence $C+VP$ is full row
rank.
\end{proof}

To conclude, the following is the algorithm to solve the control
problem with input-output partitioning constraint.

\begin{algorithm}
The following steps provide a solution to the problem if and only
if it is solvable.\newline1. Verify if the specification
$\mathcal{S}$ is regularly achievable. If so, go to step 2,
otherwise the problen is not solvable.\newline2. Construct the
canonical controller for this problem, denote it as
$\mathcal{C}_{\mathrm{{can}}}$.\newline3. Construct a regular
controller
$\mathcal{C}\in\mathfrak{C}_{\mathrm{can}}^{\mathrm{reg}}$.
Theorem \ref{t2} guarantees that this can be done. The controller
$\mathcal{C}$ and the control manifest behavior $\mathcal{P}_{c}$
can be represented in the form
of%
\begin{align}
\mathcal{C}  &  =\left\{  (u,y)~|~C_{1}\left(  \frac{d}{dt}\right)
u+C_{2}\left(  \frac{d}{dt}\right)  y=0\right\}  ,\\
\mathcal{P}_{c}  &  =\left\{  (u,y)~|~P_{1}\left(
\frac{d}{dt}\right) u+P_{2}\left(  \frac{d}{dt}\right)
y=0\right\}  .
\end{align}
\newline4. Verify if%
\begin{equation}
\text{rank}\left[
\begin{array}
[c]{c}%
M_{1}\\
P_{1}%
\end{array}
\right]  \geq\mathtt{p}(\mathcal{C}), \label{18}%
\end{equation}
where $\mathtt{p}(\mathcal{C})$ denotes the number of output
variables of $\mathcal{C}$. If (\ref{18})is satisfied, go to step
5, otherwise the problem is not solvable.\newline5. Compute a $V$
such that $C_{1}+VP_{1}$ is full row rank. The existence of such
$V$ is guaranteed by Lemma \ref{l3}. A controller
that solves the control problem is given by%
\begin{equation}
\mathcal{C}^{\prime}=\left\{  (u,y)~|~\left[
\begin{array}
[c]{cc}%
C_{1}+VP_{1} & C_{2}+VP_{2}%
\end{array}
\right]  \left(  \frac{d}{dt}\right)  \left[
\begin{array}
[c]{c}%
u\\
y
\end{array}
\right]  =0\right\}  .
\end{equation}

\end{algorithm}

\section{Concluding remarks}

We discuss a result in the field of behavioral control theory for
linear systems. The main result of the paper is a parametrization
of all regular
controllers that are equivalent to the canonical controller $\mathfrak{C}%
_{\mathrm{can}}^{\mathrm{reg}}$. This class of controllers has two
nice
properties:\newline(i) All its members are regular controllers, and\newline%
(ii) it acts as an upperbound to other regular controllers. This
means, any
regular controller is contained in an element of $\mathfrak{C}_{\mathrm{can}%
}^{\mathrm{reg}}$.

The special properties of the class $\mathfrak{C}_{\mathrm{can}}%
^{\mathrm{reg}}$ and its parametrization is used to solve two
control problems in the behavioral framework. The first control
problem is related to designing a regular controller that uses as
few control variable as possible. The second problem is about
designing a regular controller that satisfies a predefined
input-output partitioning.

The use of the parametrization of
$\mathfrak{C}_{\mathrm{can}}^{\mathrm{reg}}$ is not necessarily
limited to the above mentioned problems. An interesting problem
is, for example, to use the parametrization to construct a regular
controller with as small MacMillan degree as possible
\cite{Willems97}. Such a result can potentially lead to the
solution to the long standing problem of regular feedback
implementability \cite{Trentelman02c}.

\bibliography{adhs}
\bibliographystyle{ieeetr}

\end{document}